\documentstyle[12pt,amsmath,amssymb]{article}
\begin{document}

\newtheorem{lem}{Lemma}[section]
\newtheorem{theorem}{Theorem}[section]
\newtheorem{prop}{Proposition}[section]
\newtheorem{rem}{Remark}[section]
\newtheorem{define}{Definition}[section]
\newtheorem{cor}{Corollary}[section]
\allowdisplaybreaks

\makeatletter\@addtoreset{equation}{section}\makeatother
\def\theequation{\arabic{section}.\arabic{equation}}

\newcommand{\One}{\chi}

\newcommand{\D}{{\cal D}}
\newcommand{\N}{{\Bbb N}}
\newcommand{\C}{{\Bbb C}}
\newcommand{\Z}{{\Bbb Z}}
\newcommand{\R}{{\Bbb R}}
\newcommand{\la}{\langle}
\newcommand{\ra}{\rangle}
\newcommand{\rom}[1]{{\rm #1}}
\newcommand{\FC}{{\cal F}C_{\mathrm b}(C_0(X),\Gamma)}
\newcommand{\eps}{\varepsilon}
\newcommand{\dd}{\overset{{.}{.}}}
\newcommand{\fii}{\varphi}

\newcommand{\supp}{\operatorname{supp}}

\def\stackunder#1#2{\mathrel{\mathop{#2}\limits_{#1}}}
\newcommand{\FCo}{{\cal F}C_{\mathrm b}^\infty({\cal D},\dd\Gamma)}

\newcommand{\EG}{{\cal E}_{\mathrm G}}
\newcommand{\EK}{{\cal E}_{\mathrm K}}
\newcommand{\La}{\Lambda}
\newcommand{\ga}{\gamma} \newcommand{\Ga}{\Gamma}

\renewcommand{\author}[1]{\medskip{\large #1}\par\medskip}
\begin{center}{\Large \bf
A note on equilibrium Glauber  and Kawasaki dynamics for fermion point processes
}\end{center}

{\large Eugene Lytvynov}\\ Department of Mathematics,
University of Wales Swansea, Singleton Park, Swansea SA2 8PP, U.K.\\
e-mail: \texttt{e.lytvynov@swansea.ac.uk}\vspace{2mm}

{\large Nataliya Ohlerich}\\
 Fakult\"at f\"ur Mathematik, Universit\"at
Bielefeld, Postfach 10 01 31, D-33501 Bielefeld, Germany;
 BiBoS, Univ.\ Bielefeld,
Germany.\\ e-mail: \texttt{nturchyn@math.uni-bielefeld.de }

{\small

\begin{center}
{\bf Abstract}
\end{center}

\noindent  We construct two types of equilibrium dynamics of infinite particle systems in a
locally compact Polish space $X$, for which certain fermion point processes are invariant.
The Glauber dynamics is a birth-and-death process in $X$, while 
in the case of the Kawasaki dynamics interacting particles randomly hop over $X$.  
We establish conditions on  generators of both dynamics under which corresponding conservative Markov processes exist. 
 } \vspace{2mm}

\vspace{1.5mm}

 \vspace{1.5mm}

\section{Introduction} 

Let $X$ be a locally compact Polish space. Let $\nu$ be a Radon measure on $X$ and let $K$ be   
a  linear, Hermitian, locally trace class operator on  $L^2(X,\nu)$ for  which
 $\pmb 0\le K\le \pmb 1$. Then $K$ is an integral operator and we denote by $K(\cdot,\cdot)$ the integral kernel of $K$. 
 
 Let $\Gamma=\Gamma_X$ denote the space of all locally finite subsets (configurations) in $X$.
A fermion point process  (also called determinantal
point process) corresponding to $K$ is a probability measure on $\Gamma$
whose correlation functions are given by 
\begin{equation}
k^{(n)}_\mu(x_1,\dots,x_n)=\det(K(x_i,x_j))_{i,j=1}^n.
\end{equation}

Fermion point processes were introduced  by Macchi  \cite{Macchi} (see also Girard~\cite{Girard} and Menikoff~\cite{M2}).  These  processes 
naturally arise in quantum mechanics, statistical mechanics, random matrix theory, and   representation theory, see e.g.\ \cite{BO,ST,TI1,So2000RMS} and the references therein. 

In \cite{Spohn1}, Spohn  investigated a diffusion fynamics on the configuration space $\Gamma_{\mathbb R}$ for which the fermion process
corresponding to the Dyson (sine) kernel 
$$K(x,y)=\sin(x-y)/(x-y)$$ is an invariant measure.  

 In the case where the operator $K$ satisfies $K< \pmb 1$, Georgii and Yoo \cite{GeYoo05} (see also \cite{Yoo06}) investigated Gibbsianness of fermion point processes. In particular, they proved that every fermion process with $K$ as above possesses Papangelou (conditional) intensity. 

Using Gibbsianness of fermion point processes, Yoo
\cite{Yoo} constructed an equilibrium diffusion dynamics on the configuration space over $\R^d$, which has the fermion process as invariant measure. This Markov process is an analog of the gradient stochastic dynamics which has the standard Gibbs measure corresponding to a potential of pair interaction as invariant measure (see e.g.\ \cite{111}).

On the other hand, in the case of a standard Gibbs measure, one considers further classes of equilibrium processes on the configuration space: the so-called Glauber and Kawasaki dynamics in continuum.

The generator of the Glauber dynamics for a continuous particle system in $\R^d$, which is a birth-and-death
process,  is informally given by the formula 
\begin{equation}\label{gdtrd}(H_{\mathrm G}F)(\gamma)=\sum_{x\in\gamma}d(x,\gamma\setminus x)(D^-_xF)(\gamma)+
\int _{\R^d}b(x,\gamma) (D^+_xF)(\gamma)\,dx,\end{equation}
 where 
\begin{equation}\label{d}(D_x^-F)(\gamma)=F(\gamma\setminus x)-F(\gamma),\quad 
 (D_x^+F)(\gamma)=F(\gamma\cup x)-F(\gamma).\end{equation}
Here and below, for simplicity of notations, we just write $x$ instead of $\{x\}$.
The coefficient $d(x,\gamma\setminus x)$ describes the rate at which the particle $x$ of the configuration $\gamma$ dies, while $b(x,\gamma)$ describes the rate at which, given the configuration $\gamma$, a new
particle is born  at $x$.     

The Kawasaki dynamics of continuous particles is a process in which particles randomly hop
over the space $\R^d$.  
 The generator of such a process is then informally given by 
\begin{equation}\label{serserser}(H_{\mathrm K}F)(\gamma)=\sum_{x\in\gamma}\int_{\R^d} c(x,y,\gamma\setminus x)(D^{-+}_{xy}F)(\gamma)\,dy,\end{equation} where 
\begin{equation}\label{dd}
(D_{xy}^{-+}F)(\gamma)=F(\gamma\setminus x\cup y)-F(\gamma)
\end{equation}
and the coefficient $c(x,y,\gamma\setminus x)$ describes the rate at which the particle $x$ of the configuration $\gamma$ 
jumps to $y$. 

Glauber and Kawasaki dynamics of continuous particle systems in infinite volume which have a standard Gibbs measure as symmetrizing measure were
constructed in \cite{KL,KLRprep}, see also \cite{BCC,G1,G2,HS,P,Wu}. For further studies on Glauber and Kawasaki dynamics, we refer to \cite{FKL,Grothaus,KKL,KMZh}. 

So, the aim of the present note is to prove the existence of Glauber and Kawasaki dynamics of a continuous particle system which has a fermion point process as invariant measure.
Our choice of dynamics seems to be absolutely natural for a fermion system. 
We recall that Shirai and Yoo \cite{ShYoo} already constructed  a Glauber dynamics on the lattice which has a fermion point process on the lattice as invariant measure.  

Using the theory of Dirichlet forms (see e.g.\ \cite{MR}), we will construct conservative Markov processes on $\Gamma$ with {\it cadlag} paths which have a fermion measure $\mu$ as symmetrizing, hence invariant  measure. Furthermore, we will  discuss the explicit form of the $L^2(\mu)$-generators of these process on a set of 
 cylinder functions. These generators will have the form  \eqref{gdtrd} in the case of Glauber dynamics, and \eqref{serserser} in the case of Kawasaki dynamics
 (with $\R^d$ replaced by a general topological space $X$). Since we essentially use the Papangelou intensity of the fermion point process, our study here is restricted by the assumption that $K<\pmb 1$.

Throughout the paper, we formulate our results for both dynamics, and give the proofs only in the case
of the Kawasaki dynamics. The reason is that the proofs in the Glauber case  are quite similar to, and simpler than the corresponding proofs for the Kawasaki dynamics. 

Finally, let us mention some open problems, which will be topics of our further research:

\begin{enumerate}
\item Construction of Glauber and Kawasaki dynamics for fermion point processes  in the case where  1 belongs to the spectrum of the  operator $K$
(which is, e.g., the case when $K$ has Dyson kernel).

\item Finding a core  for the generator of the dynamics (compare with \cite{KL}).

\item Deciding whether the generator of the Glauber dynamics for a fermion point process has a spectral gap. 
(Recall that, in the case of a standard Gibbs measure with a positive  potential of pair interaction, there is a Glauber dynamics whose generator has a spectral gap \cite{KL}).

\item Study of different types of scalings of Glauber and Kawasaki dynamics, in particular, difffusion approximation for Kawasaki dynamics (compare with \cite{FKL,Grothaus,KKL}).
\end{enumerate}

\section{Fermion (determinantal) point processes}\label{fsc}

Let $X$ be a locally compact, second countable Hausdorff
topological space. Recall that such a space is known to be Polish. We denote by ${\mathcal B}(X)$ the Borel $\sigma$-algebra
in $X$, and by ${\mathcal B}_{0}(X)$ the collection of all sets from
${\mathcal B}(X)$ which are relatively compact. 
We fix a Radon, non-atomic measure $\nu$
on  $(X, {\cal B}(X))$.

The configuration space $\Gamma:=\Gamma_{X}$ over $X$ is defined as the set of all subsets of $X$ which are
locally finite: $$\Gamma:=\big\{\gamma\subset X : \,
|\gamma_\Lambda|<\infty\text{ for each  }\Lambda\in{\mathcal B}_0(X )\big\},$$ where $|\cdot|$ denotes the cardinality of a set
and $\gamma_\Lambda:= \gamma\cap\Lambda$. One can identify any
$\gamma\in\Gamma$ with the positive Radon measure
$\sum_{x\in\gamma}\eps_x\in{\cal M}(X)$,  where  $\eps_x$ is
the Dirac measure with mass at $x$,
$\sum_{x\in\varnothing}\varepsilon_x{:=}$zero measure, and ${\cal
M}(X)$
 stands for the set of all
positive  Radon  measures on   ${\cal
B}(X)$. The space $\Gamma$ can be endowed with the relative
topology as a subset of the space ${\cal M}(X)$ with the vague
topology, i.e., the weakest topology on $\Gamma$ with respect to
which  all maps $$\Gamma\ni\gamma\mapsto\la
f,\gamma\ra:=\int_{X} f(x)\,\gamma(dx) =\sum_{x\in\gamma}f(x),\qquad
f\in C_0(X),$$ are continuous. Here, $C_0(X)$ is the space
of all continuous,  real-valued functions on $X$ with compact
support. We will  denote  the Borel
$\sigma$-algebra on $\Gamma$ by ${\cal B}(\Gamma)$. A point process $\mu$ is a 
probability measure on $(\Gamma,{\cal B}(\Gamma))$.

A point process $\mu$  is said to have correlation functions
if, for any $n\in\N$, there exists a non-negative, measurable,
symmetric function $k_\mu^{(n)}$ on $X^n$  such
that, for any measurable, symmetric function
$f^{(n)}:X^n\to[0,+\infty]$,
\begin{align*} &\int_\Gamma \sum_{\{x_1,\dots,x_n\}\subset\gamma}
f^{(n)}(x_1,\dots,x_n)\,\mu(d\gamma)\notag\\&\qquad =\frac1{n!}\,
\int_{X^n} f^{(n)}(x_1,\dots,x_n)
k_\mu^{(n)}(x_1,\dots,x_n)\,\nu(dx_1)\dotsm
\nu(dx_n).
\end{align*}

Let $K$ be a linear, bounded, Hermitian  operator on the space $L^2(X,\nu)$ (real or complex) which satisfies the following assumptions:

\begin{enumerate}
\item $K$ is locally of trace class, i.e., 
$$
\operatorname{Tr}(P_\Lambda K P_\Lambda)<\infty \quad \text{for all } \Lambda\in{\cal B}_{0}(X),
$$
where $P_\Lambda$ denotes the  operator of multiplication by the indicator function $\One_\Lambda$ of the set $\Lambda$.

\item We have $\pmb 0\le K\le \pmb 1$. 
\end{enumerate}

Under the above assumptions,  $K$ is an integral operator, and its kernel can be chosen as $$K(x,y)=\int_X
K_1(x,z)K_1(z,y)\,\nu(dz),$$ where $K_1(\cdot,\cdot)$ is any version of the kernel of the integral operator $\sqrt K$, \cite{LM}
(see also \cite[Lemma A.4]{GeYoo05}).

A point process $\mu$ having correlation functions 
$$
k^{(n)}_\mu(x_1,\dots,x_n)=\det(K(x_i,x_j))_{i,j=1}^n.
$$
is called the fermion (or determinantal) point process corresponding to the operator  $K$.
Under the above assumptions on $K$, such a  point process $\mu$ exists and is unique, see
e.g.\  \cite{LM,Macchi,So2000RMS,ST}. 

Using the definition of a fermion process, it can be easily checked that $\mu$ has all local moments finite, i.e., 
\begin{equation*}
\int_\Gamma \la f,\gamma\ra^n\,\mu(d\gamma)<\infty,\qquad
f\in C_0(X),\ f\ge0,\  n\in\N.\end{equation*} 

In what follows,  we will always assume that the operator $K$ is strictly less than
$\pmb 1$, i.e., 1 does not belong to the spectrum of $K$. Then, as has been shown  in 
\cite{GeYoo05}, the fermion process $\mu$ has Papangelou (conditional) intensity. That is, there exists a measurable
function $r:X\times\Gamma\to[0,+\infty]$ such that
 \begin{equation}
\int_\Gamma \mu(d\gamma)\int_{X} \gamma(dx)
F(x,\gamma) =\int_\Gamma \mu(d\gamma)\int_{X}
\nu(dx)\, r(x,\gamma) F(x,\gamma\cup x)\label{mecke}
\end{equation}
for any measurable function $F:X\times \Gamma \to[0,+\infty].$

\begin{rem}\rom{Let us briefly explain the construction of the Papangelou intensity $r(x,\gamma)$, following
\cite{GeYoo05}. 

For each $\Lambda\in {\mathcal B}_0(X)$, consider $K_\Lambda:=P_\Lambda K P_\Lambda$ as an operator in 
$L^2(\Lambda,\nu)$ and let $J_{[\Lambda]}:=K_\Lambda (\pmb 1-K_\Lambda)^{-1}$. Denote 
by $J_{[\Lambda]}(\cdot,\cdot)$ the kernel of the operator $J_{[\Lambda]}$ (chosen analogously to the kernel of $K$). For any
$\gamma\in\Gamma$, set $$ \operatorname{det}J_{[\Lambda]}(\gamma_\Lambda,\gamma_\Lambda):=
\operatorname{det} \big[J_{[\Lambda]}(x_i,x_j)\big]_{i,j=1}^m, $$ with $\gamma_\Lambda=\{x_1,\dots,x_m\}$
being any numeration of points of $\gamma_\Lambda$ (in the case $\gamma_\Lambda=\varnothing$, set  $\operatorname{det}J_{[\Lambda]}(\varnothing,\varnothing):=0$). Now, for any $x\in\Lambda$ and $\gamma\in\Gamma$,
set $$ r_\Lambda(x,\gamma_\Lambda):=\frac{\operatorname{det}J_{[\Lambda]}(x\cup\gamma_\Lambda,x\cup\gamma_\Lambda)}
{\operatorname{det}J_{[\Lambda]}(\gamma_\Lambda,\gamma_\Lambda)}\,, $$
where the expression on the right hand side is assumed to be zero if $\operatorname{det}J_{[\Lambda]}(\gamma_\Lambda,\gamma_\Lambda)=0$.  

Let $\{\Lambda_n\}_{n\in\mathbb N}$ be any sequence in ${\mathcal B}_0(X)$ that increases to $X$.
Then $r(x,\gamma)$ is a $\nu\otimes\mu$-a.e.\ limit of $r_{\Lambda_n}(x,\gamma_{\Lambda_n})$ as $n\to\infty$. 
}\end{rem}

Set $J:=K(\pmb 1-K)^{-1}$. The operator $J$
is integral and we choose its kernel $J(\cdot,\cdot)$ analogously to choosing the kernel of $K$. 
 Note that 
$$ \operatorname{Tr}(P_\Lambda J P_\Lambda)=\int_\Lambda J(x,x)\,\nu(dx)<\infty.$$

The following proposition 
is a direct corollary of Theorem 3.6 and   \cite[Lemma~A.1]{GeYoo05}.

\begin{prop}\label{gtyfy}
We have, for $\nu\otimes\mu$-a.e.\ $(x,\gamma)\in X\times\Gamma$:
$$ r(x,\gamma)\le J(x,x).$$
\end{prop}

\section{Equilibrium dynamics}\label{rfyje}

In what follows, we will consider a fermion point process $\mu$ corresponding to an operator $K$
 as defined in Section~\ref{fsc}. We introduce the set $\FC$
of all functions of the form $$\Gamma\ni \gamma\mapsto F(\gamma)=g(\la\varphi_1,\gamma\ra,\dots , \la\varphi_N,\gamma\ra), $$ where $N\in\N$, $\varphi_1,\dots,\varphi_N\in C_0(X)$ and $g\in C_{\mathrm b}(\R^N)$. Here,  $C_{\mathrm b}(\R^N)$ denotes the set of all continuous, bounded functions on $\R^N$. 

For a function $F:\Gamma\to\R$, $\gamma\in\Gamma$, $x,y\in X$, we introduce the  notations $(D_x^-F)(\gamma)$, 
$(D_x^+F)(\gamma)$, and  
$(D_{xy}^{-+}F)(\gamma)$ by \eqref{d} and  \eqref{dd}, respectively. We consider measurable mappings \begin{gather*} X\times\Gamma \ni (x,\gamma)\mapsto d(x,\gamma)\in[0,\infty),\\
X\times X\times\Gamma\ni (x,y,\gamma)\mapsto c(x,y,\gamma)\in[0,\infty).\end{gather*}
Assume that 
\begin{equation}\label{condC}
c(x,y,\gamma)=c(x,y,\gamma)\One_{\{r(x,\gamma)>0,\, 
r(y,\gamma)>0\}},\qquad x,y\in X,\ \gamma\in\Gamma.
\end{equation}

\begin{rem}\rom{ As we will see below, the coefficient $c(x,y,\gamma\setminus x)$ describes the rate of the jump of  particle $x\in\gamma$ to $y$. If $r(y,\gamma\setminus x)=0$, then the relative energy of interaction between  the configuration $\ga\setminus x$ and point $y$ is $+\infty$, so that the particle $x$ cannot jump to $y$, i.e., $c(x,y,\gamma\setminus x)$ should be equal to zero. A symmetry reason also implies that we should have $c(x,y,\gamma\setminus x)=0$ if $r(x,\gamma\setminus x)=0$, i.e., if the relative energy of interaction between $x\in\gamma$ and the rest of configuration is $\infty$. 

}\end{rem}

Further, we assume that, for each  $\Lambda\in{\cal B}_0( X)$,
\begin{gather}
\int_\Gamma\mu(d\gamma)\int_\Lambda\gamma(dx)d(x,\gamma\setminus x)<\infty,\label{gyiu}\\
\int_\Gamma\mu(d\gamma)\int_X\gamma(dx)\int_X \nu(dy)c(x,y,\gamma\setminus x)(\One_\Lambda(x)+\One_\Lambda(y))<\infty.
\label{gyu}\end{gather}


We define bilinear forms
\begin{align} \EG(F,G):=& \int_\Gamma\mu(d\gamma)\int_X \gamma(dx)\,d(x,\gamma\setminus x) 
(D_x^-F)(\gamma) (D_x^-G)(\gamma),\label{1}\\
\EK (F,G):=& \int_\Gamma \mu(d\gamma)\int_X \gamma(dx)\int_X \nu(dy)\, c(x,y,\gamma\setminus x)
 (D_{xy}^{-+}F)(\gamma) (D_{xy}^{-+}G)(\gamma),\label{2}
\end{align}
where $F,G\in\FC$. 

We note that, for any $F\in\FC$, there exist  $\Lambda\in{\cal B}_0( X)$ and $C>0$ such that 
$$ |(D^-_xF)(\gamma)|\le C\One_\Lambda(x),\quad 
|(D_{xy}^{-+}F)(\gamma)|\le C (\One_\Lambda(x)+\One_\Lambda(y)),\qquad \gamma\in\Gamma,\ x,y\in X.$$
Therefore, by assumptions \eqref{gyiu}, \eqref{gyu}
the right-hand sides of formulas \eqref{1} and \eqref{2} are well-defined and finite.

Using \eqref{mecke} and \eqref{condC}, we have, for any $F\in\FC$:
\begin{align*}
\EK(F)&=\int_\Gamma\mu(d\gamma)\int_X\nu(dx)\int_X\nu(dy)\,r(x,\gamma)c(x,y,\gamma)\\
&\quad\times\One_{\{r(y,\gamma)>0)\}}\frac{r(y,\gamma)}{r(y,\gamma)}(F(\gamma\cup
y)-F(\gamma\cup x))^2\\
&=\int_\Gamma \mu(d\gamma)\int_X\nu(dx)\int_X \gamma(dy)\,
r(x,\gamma\setminus y)c(x,y,\gamma\setminus y)\\
&\quad\times\One_{\{r(y,\gamma\setminus y)>0\}}\frac1{r(y,\gamma\setminus y)}(D_{yx}^{-+}F)^2(\gamma).
\end{align*}
Here, we used the notation $\EK(F):=\EK(F,F)$.
Therefore, for any $F,G\in\FC$,
$$ \EK(F,G)=\int_\Gamma\mu(d\gamma)\int_X\gamma(dx)\int_X
\nu(dy)\tilde c(x,y,\gamma\setminus x)(D_{xy}^{-+}F)(\gamma)(D_{xy}^{-+}G)(\gamma),$$
where
$$\tilde c(x,y,\gamma):=\frac12\bigg(
c(x,y,\gamma)+c(y,x,\gamma)\One_{\{r(x,\gamma)>0\}}\,\frac{r(y,\gamma)}{r(x,\gamma)}
\bigg). $$
Note that, by \eqref{condC}, we have $\tilde{\tilde c}(x,y,\gamma)=\tilde c(x,y,\gamma)$. Therefore, without loss of generality, in what follows we will assume that $\tilde c(x,y,\gamma)=c(x,y,\gamma)$, i.e., \begin{equation}\label{hyhgfuy}
r(x,\gamma)c(x,y,\gamma)=r(y,\gamma)c(y,x,\gamma).\end{equation}

\begin{lem} We have ${\cal E}_\sharp(F,G)=0$ for all $F,G\in\FC$ such that $F=0$ $\mu$-a\rom.e\rom.\rom, $\sharp={\mathrm G}, {\mathrm K}$\rom. 
\end{lem}

\noindent {\it Proof.}  It suffices to show that, for $F\in\FC$ such that $F=0$ $\mu$-a.e., we have  
$(D_{x,y}^{-+}F)(\gamma)=0$ $\tilde\mu$-a.e.
Here, $\tilde\mu$ is the measure on $X\times X\times\Gamma$ defined by \begin{equation}\label{yguuyuy}
\tilde\mu(dx,dy,d\gamma):=c(x,y,\gamma\setminus x)
\gamma(dx)\,\nu(dy)\,\mu(d\gamma).\end{equation}

For any $F$ as above, we evidently have that $F(\gamma)=0$ $\tilde\mu$-a.e. 
 Next, by \eqref{mecke} and \eqref{condC}
\begin{align}\label{yuuy} &\int_\Gamma \mu(d\gamma)\int_\Lambda\gamma(dx)
\int_\Lambda \nu(dy)|F(\gamma\setminus x\cup y)|c(x,y,\gamma\setminus x)\notag\\
&\qquad=\int_\Gamma \mu(d\gamma)\int_\Lambda\nu(dx)
\int_\Lambda \nu(dy)r(x,\ga)|F(\gamma\cup y)|c(x,y,\gamma)\chi_{\{r(y,\ga)>0\}}
\frac{r(y,\ga)}{r(y,\ga)}\notag\\
&\qquad=\int_\Gamma \mu(d\gamma)\int_\Lambda\nu(dx)
\int_\Lambda \ga(dy)|F(\gamma)|c(x,y,\gamma\setminus y)\,\frac {r(x,\ga\setminus y)}{r(y,\gamma\setminus y)}
\,\chi_{\{r(y,\ga\setminus y)>0\}}\notag\\
&\qquad=\int_\Gamma\mu(d\gamma) \int_\Lambda\nu(dx)\int_\Lambda \ga(dy)|F(\gamma)|c(x,y,\gamma\setminus y)
\frac {r(x,\ga\setminus y)}{r(y,\gamma\setminus y)}.\end{align}
Since $F$ is bounded,   by \eqref{gyu} the integral in \eqref{yuuy} is finite. Therefore, 
\begin{equation}\label{ghiuygiyg}|F(\gamma)|\,
\frac {r(x,\ga\setminus y)}{r(y,\gamma\setminus y)}
<\infty \quad \text{for $\tilde \mu$-a.a.\ $(x,y,\gamma)\in X\times X\times\Gamma$}.\end{equation} 
Because $F=0$ $\tilde\mu$-a.e., by \eqref{yuuy} and \eqref{ghiuygiyg}, 
$F(\gamma\setminus x\cup y)=0$ $\tilde\mu$-a.e.
\quad $\square$\vspace{2mm}

\begin{lem}\label{uiwfg}
\rom{1)} The bilinear form $(\EG,\FC)$ is closable on $L^2(\Gamma,\mu)$	
 and its closure will be denoted by $(\EG,D(\EG))$\rom.

\rom{2)} Assume that, for some $u\in\mathbb R$,
\begin{equation} \label{cond2}\int_\Lambda\nu(dx)\int_\Lambda\gamma(dy)r(x,\gamma\setminus y) r
(y,\gamma\setminus y)^{u}\One_{\{r(y,\gamma\setminus y)>0\}}c(x,y,\gamma\setminus y)\in L^2(\Gamma,\mu)\end{equation}
for all $\Lambda\in{\mathcal B}_0(X)$.
Then the bilinear form $(\EK,\FC)$ is closable on $L^2(\Gamma,\mu)$	 and its closure will be denoted by $(\EK,D(\EK))$\rom. 
\end{lem} 

\noindent {\it Proof.} Let $(F_n)_{n=1}^\infty$ be a sequence in $\FC$ such that 
$\|F_n\|_{L^2(\Gamma,\mu)}\to0$ as $n\to\infty$ and
\begin{equation}\label{rdd}\EK(F_n-F_k)\to0\quad \text{as $n,k\to\infty$}.\end{equation} 
  To prove the closability of $\EK$, it suffices to show that there exists a subsequence $\{F_{n_k}\}_{k=1}^\infty$ such that
$\EK(F_{n_k})\to0$ as $k\to\infty$.

Since $\|F_n\|_{L^2(\Gamma,\mu)}\to0$ as $n\to\infty$, there exists a subsequence $(F_n^{(1)})_{n=1}^\infty $ of  $(F_n)_{n=1}^\infty$ such that 
$F_n^{(1)}(\gamma)\to0$ for $\tilde\mu$-a.a.\ $(x,y,\gamma)\in X\times X\times \Gamma$.
Next, by \eqref{cond2}, we have, for any  $\Lambda\in{\cal B}_0(X)$,   \begin{align*}
&\int_\Gamma \mu(d\gamma)\int_\Lambda \gamma(dx)\int_\Lambda \nu(dy)
c(x,y,\gamma\setminus x)r(y,\gamma\setminus x)^{u+1}\One_{\{ r(y,\gamma\setminus x)>0\}}|F^{(1)}_n(\gamma\setminus x\cup y)|\\
&\qquad=\int_\Gamma \mu(d\gamma)\int_\Lambda\nu(dx)\int_\Lambda
\nu(dy)r(x,\gamma)c(x,y,\gamma)r(y,\gamma)^{u+1} \One_{\{ r(y,\gamma)>0\}}
|F^{(1)}_n(\gamma\cup y)|\\
&\qquad=\int_\Gamma\mu(d\gamma)\int_\Lambda\nu(dx)\int_\Lambda \gamma(dy)r(x,\gamma\setminus y)r(y,\gamma\setminus y)^{u}\\
&\quad\qquad\times\One_{\{
r(y,\gamma\setminus y)>0\}}c(x,y,\gamma\setminus y)|F^{(1)}_n(\gamma)|\\
&\qquad\le \bigg(\int_\Gamma\mu(d\gamma) |F^{(1)}_n(\gamma)|^2\bigg)^{1/2}
\bigg(\int_{\Gamma}\mu(d\gamma)\bigg(\int_\Lambda\nu(dx)\int_\Lambda \gamma(dy)r(x,\gamma\setminus y)\\
&\qquad\quad\times r(y,\gamma\setminus y)^{u}
\One_{\{r(y,\gamma\setminus y)>0\}}c(x,y,\gamma\setminus y)\bigg)^2\bigg)^{1/2}\to0\quad\text{as }n\to\infty.
\end{align*}
Therefore, there exists a subsequence $(F^{(2)}_n)_{n=1}^\infty$ of 
$(F^{(1)}_n)_{n=1}^\infty$ such that $F_n^{(2)}(\gamma\setminus  x\cup y)\to 0$ as $n\to\infty$ for $$c(x,y,\gamma\setminus x)r(y,\gamma\setminus x)^u\One_{\{r(y,\gamma\setminus x)>0\}}\gamma(dx)\nu(dy)\mu(d\gamma)\text{-a.e.\ }(x,y,\gamma)\in X\times X\times\Gamma.$$ By \eqref{condC}, the latter measure is equivalent
to $\tilde\mu$, and therefore \begin{equation}\label{igiuygtyi}
(D^{-+}_{x,y}F^{(2)}_n)(\gamma)\to0 \quad \text{for $\tilde\mu$-a.e.\ $(x,y,\gamma)\in X\times X\times\Gamma$}.\end{equation}

Now, by \eqref{igiuygtyi} and Fatou's lemma
\begin{align*}
\EK(F_n^{(2)})&=\int (D_{xy}^{-+}F_n^{(2)})(\gamma)^2\,\tilde\mu(dx,dy,d\gamma)\\&=
\int \left((D_{xy}^{-+}F_n^{(2)})(\gamma)-
\lim_{m\to\infty}
(D_{xy}^{-+}F_m^{(2)})(\gamma)\right)^2\,\tilde\mu(dx,dy,d\gamma)\\&\le 
\liminf_{m\to\infty}\int ((D_{xy}^{-+}F^{(2)}_n)(\gamma)-(D_{xy}^{-+}F^{(2)}_m)(\gamma))^2\,\tilde\mu(dx,dy,d\gamma)\\&=\liminf_{m\to\infty}\EK(F_n^{(2)}-F_m^{(2)}),
\end{align*}
which by \eqref{rdd} can be made arbitrarily small for $n$ large enough.\quad $\square$\vspace{2mm}

We will now need  the bigger  space $\dd\Gamma$ consisting of all
$\Z_+\cup\{\infty\}$-valued Radon measures on $X$ (which is Polish, see e.g.\
\cite{Ka75}). Since $\Gamma\subset\dd\Gamma$ and ${\cal
B}(\dd\Gamma)\cap\Gamma={\cal B}(\Gamma)$, we can consider $\mu $
as a measure on $(\dd\Gamma,{\cal B}(\dd\Gamma))$ and
correspondingly $({\cal E}, D({\cal E}))$ as a bilinear form on
$L^2(\dd\Gamma,\mu)$.

For the notion of a quasi-regular Dirichlet form,  appearing in
the following lemma, we refer to \cite[Chap.~I, Sect.~4 and Chap.~IV, Sect.~3]{MR}.

\begin{lem}\label{scfxgus} Under the assumption of Lemma~\rom{\ref{uiwfg}},
$(\EG,D(\EG))$ and $(\EK,D(\EK))$ are quasi-regular Dirichlet forms on
$L^2(\dd\Gamma,\mu)$.
\end{lem}

The proof of Lemma \ref{scfxgus}  is analogous to that  of \cite[Lemmas 3.3 and 3.4]{KLRprep}, so we omit it.

For the notion of an exceptional set, appearing in
the next proposition, we refer  e.g.\  to \cite[Chap.~III,
Sect.~2]{MR}.

\begin{lem}\label{fdj} The set  $\dd\Gamma\setminus\Gamma$ is 
exceptional for both $\EG$ and $\EK$\rom.
\end{lem}

\noindent {\it Proof}.  We fix any metric on $X$ which generates the topology on $X$.
For any $a\in X$ and $r>0$, we denote by $B(a,r)$ the closed ball in $X$, with center at $x$ and radius $r$.
It suffices to prove the lemma locally, i.e., to show that, for any fixed $a\in X$,
there exists $r>0$ such that $$N_a:=\{\gamma\in\ddot\Gamma: \sup_{x\in B(a,r)}\gamma(\{x\})\ge2\}$$ is $\EK$-exceptional. 
So, we fix any $a\in X$ and choose $r>0$ so that $B(a,2r)\in{\mathcal B}_0(X)$. 

By
\cite[Lemma~1]{RS98}, we need to prove that there exists a
sequence $u_n\in D(\EK)$, $n\in\N$, such that each $u_n$ is a
continuous  function on $\dd\Gamma$, $u_n\to 
\chi_{N_a}$ pointwise
as $n\to\infty$, and $\sup_{n\in\N}\EK (u_n)<\infty$.

 Fix any $n\in\mathbb N$ such that \begin{equation}\label{jhfytfty}2/n<r.
 \end{equation}
  Let $$\{B(a_k,1/n)\mid k=1,\dots,K_n\},$$ with $a_k\in B(a,r)$, $k=1,\dots,K_n$, be a finite covering of $B(a,r)$.  Let $f:\mathbb R\to\mathbb R$ be given by $f(u):=(1-|u|)\vee 0$. 
    
  For each $k=1,\dots,K_n$, we define a continuous function $f_k^{(n)}:X\to\mathbb R$ by 
$$ f_k^{(n)}(x):=f\big(n\,\operatorname{dist}(x,B(a_k,1/n))\big),\quad x\in X.$$
Here, $\operatorname{dist}(x,B)$ denotes the distance from a point $x\in X$ to a set $B\subset X$. We evidently have:
\begin{equation}\label{hyturfiuuu} \One_{B(a_k,1/n)}\le f_k^{(n)}\le \One_{B(a_k,2/n)}. \end{equation}

Let $\psi\in C_{\mathrm b}^1 ({\Bbb R})$ be such that $\One
_{[2,\infty)}\le\psi\le\One_{[1,\infty)}$ and \begin{equation}\label{yufutf} 
0\le\psi'\le 2\,
\One_{ (1,\infty)}.\end{equation}
We define a continuous function $$\dd\Gamma\ni\gamma\mapsto
u_n(\gamma){:=}\psi \left( \sup_{k\in \{1,\dots,K_n\}}\la
f_k^{(n)},\gamma\ra\right),$$ whose restriction to
$\Gamma$ belongs to $\FC$. Evidently, $u_n\to \One_{N_a}$ pointwise as
$n\to\infty$.

By \eqref{jhfytfty}, \eqref{hyturfiuuu}, \eqref{yufutf}, and  the mean value theorem, we have, for each $\gamma\in\Gamma$, $x\in\gamma$, $y\in X\setminus 
\gamma$, 
\begin{align}
(D_{xy}^{-+}u_n)^2(\gamma)
&\le 4\bigg(\sup_{k\in\{1,\dots,K_n\}}\la f_k^{(n)},\gamma\setminus x\cup y\ra- \sup_{k\in\{1,\dots,K_n\}}\la f_k^{(n)},\gamma\ra\bigg)^2\notag\\
&\le 4
\sup_{k\in\{1,\dots,K_n\}} |\la f_k^{(n)},\gamma\setminus x\cup y\ra- \la f_k^{(n)},\gamma\ra|^2\notag\\
&\le 8\bigg(\sup_{k\in\{1,\dots,K_n\}}f_k^{(n)}(x)^2+\sup_{k\in\{1,\dots,K_n\}}f_k^{(n)}(y)^2\bigg)\notag\\ 
&\le 8 \bigg(\sup_{k\in\{1,\dots,K_n\}}
\One_{B(a_k,2/n)}(x)
+\sup_{k\in\{1,\dots,K_n\}}\One_{B(a_k,2/n)}(y)\bigg)\notag\\
&\le 8(\One_{B(a,2r)}(x)+\One_{B(a,2r)}(y) ).\notag
\end{align} 
Hence, by \eqref{gyu}, $$\sup_{n}\EK(u_n)<\infty, $$
which implies the lemma.\quad $\square$

We now have the main result of this paper.

\begin{theorem}\label{8435476} 
Let \eqref{gyiu}, respectively   \eqref{gyu} and  \eqref{cond2},  hold. Let $\sharp={\mathrm G},{\mathrm K}$\rom. Then we have\rom:
 
\rom{1)} 
There exists a conservative Hunt process 
$${\bf M}^\sharp=({\pmb{ \Omega}}^\sharp,{\bf F}^\sharp,({\bf F}^\sharp_t)_{t\ge0},({\pmb
\Theta}^\sharp_t)_{t\ge0}, ({\bf X}^\sharp(t))_{t\ge 0},({\bf P
}^\sharp_\gamma)_{\gamma\in\Gamma})$$ on $\Gamma$ \rom(see e\rom.g\rom.\
\rom{\cite[p.~92]{MR})} which is properly associated with $({\cal
E}_\sharp,D({\cal E}_\sharp))$\rom, i\rom.e\rom{.,} for all \rom($\mu$-versions
of\/\rom) $F\in L^2(\Gamma,\mu)$ and all $t>0$ the function
\begin{equation}\label{zrd9665} \Gamma\ni\gamma\mapsto
p^\sharp_tF(\gamma){:=}\int_{\pmb\Omega} F({\bf X}^\sharp(t))\, d{\bf
P}^\sharp_\gamma\end{equation} is an ${\cal E}_\sharp$-quasi-continuous version
of $\exp(-t{H}_\sharp)F$\rom, where $(H_\sharp,D(H_\sharp))$ is the generator of $({\cal
E}_\sharp,D({\cal E}_\sharp))$\rom.  $\bf M^\sharp$ is up to $\mu$-equivalence unique
\rom(cf\rom.\ \rom{\cite[Chap.~IV, Sect.~6]{MR}).} In
particular\rom, ${\bf M}^\sharp$ is $\mu$-symmetric \rom(i\rom.e\rom{.,}
$\int G\, p^\sharp_tF\, d\mu=\int F \, p^\sharp_t G\, d\mu$ for all
$F,G:\Gamma\to{\Bbb R}_+$\rom, ${\cal B}(\Gamma)$-measurable\rom)\rom, so
 has $\mu$ as an invariant measure\rom.

\rom{2)} ${\bf M}^\sharp$ from \rom{1)} is  up to $\mu$-equivalence
\rom(cf\rom.\ \rom{\cite[Definition~6.3]{MR}}\rom) unique between
all Hunt processes ${\bf M}'=({\pmb{ \Omega}}',{\bf F}',({\bf
F}'_t)_{t\ge0},({\pmb \Theta}'_t)_{t\ge0}, ({\bf X}'(t))_{t\ge
0},({\bf P }'_\gamma)_{\gamma\in\Gamma})$ on $\Gamma$ having $\mu$
as  invariant measure and solving the martingale problem for
$(-H_\sharp, D(H_\sharp))$\rom, i\rom.e\rom.\rom, for all $G\in D(H_\sharp)$
$$\widetilde G({\bf X}'(t))-\widetilde G({\bf X}'(0))+\int_0^t (H_\sharp
G)({\bf X}'(s))\,ds,\qquad t\ge0,$$ is an $({\bf
F}_t')$-martingale under ${\bf P}_\gamma'$ for ${\cal
E}_\sharp$-q\rom.e\rom.\ $\gamma\in\Gamma$\rom. \rom(Here\rom,
$\widetilde G$ denotes an ${\cal E}_\sharp$-quasi-continuous version of $G$\rom,
cf\rom. \rom{\cite[Ch.~IV, Proposition~3.3]{MR}.)}
\end{theorem}

\begin{rem}\rom{
 In Theorem~\ref{8435476}, ${\bf M}^\sharp$ can be taken canonical, i.e., $\pmb\Omega^\sharp$ is
the set of all {\it cadlag} functions $\omega:[0,\infty)\to
\Gamma$ (i.e., $\omega$ is right continuous on $[0,\infty)$ and
has left limits on $(0,\infty)$), ${\bf
X}^\sharp(t)(\omega){:=}\omega(t)$, $t\ge 0$, $\omega\in\pmb\Omega^\sharp$,
$({\bf F}^\sharp_t)_{t\ge 0}$ together with ${\bf F}^\sharp$ is the corresponding
minimum completed admissible family (cf.\
\cite[Section~4.1]{Fu80}) and ${\pmb \Theta}^\sharp_t$, $t\ge0$, are the
corresponding natural time shifts.
 }\end{rem}

\noindent {\it Proof of Theorem\/}~\ref{8435476}. The first part
of the theorem follows from
 Lemmas~\ref{scfxgus}, \ref{fdj}, the fact that $1\in D({\cal E}_\sharp)$ and ${\cal E}_\sharp(1,1)=0$, $\sharp={\mathrm G}, {\mathrm K}$,  and
\cite[Chap.~IV, Theorem~3.5 and Chap.~V, Proposition~2.15]{MR}.
The second part follows directly from the proof of
\cite[Theorem~3.5]{AR}.\quad $\square$
\vspace{2mm}

Now we will derive explicit formulas for the generators of $\EG$ and $\EK.$ However, for this, we will demand stronger conditions on the coefficients $d(x,\gamma)$ and $c(x,y,\gamma)$.

\begin{theorem} \label{gener} 
\rom{1)} Assume that\rom, for each  $\Lambda\in{\cal B}_0( X)$\rom, \begin{align}
\int_\Lambda \gamma(dx) d(x,\gamma\setminus x)&\in L^2(\Gamma,\mu),\notag\\ 
\int_\Lambda \nu(dx)b(x,\gamma)&\in L^2(\Gamma,\mu)\label{igyig},\end{align}
where \begin{equation}\label{oii} b(x,\gamma):= r(x,\gamma)d(x,\gamma),
\qquad x\in X, \gamma\in\Gamma.\end{equation} 
Then, for each $F\in\FC$,\begin{equation}\label{genG}
(H_{\mathrm G}F)(\gamma)=-\int_{X}
\nu(dx)\, b(x,\ga)(D^+_xF)(\gamma)-
\int_{X}\gamma(dx)\,d(x,\ga\setminus x)(D^-_xF)(\gamma)\qquad \text{\rom{$\mu$-a.e.}}
\end{equation}
and $(H_{\mathrm G},D(H_{\mathrm G}))$ is the  Friedrichs extension of  $(H_{\mathrm G},\FC)$
in $L^2(\Gamma,\mu)$.

\rom{2)}   Assume that\rom, for each  $\Lambda\in {\cal B}_0( X)$\rom, \begin{equation}\label{jkg}
\int_X\gamma(dx)\int _X \nu(dy)\, c(x,y,\gamma\setminus x)(\One_\Lambda(x)+\One_\Lambda(y))\in L^2(\Gamma,\mu). 
\end{equation}
Then, for each $F\in\FC$,
\begin{equation}\label{genK}
(H_{\mathrm K}F)(\gamma)=-2\int_{X}\gamma(dx)\int _X \nu(dy)  
c(x, y,\gamma\setminus x) (D^{-+}_{xy}F)(\gamma)\qquad \text{\rom{$\mu$-a.e.}}
\end{equation}
and $(H_{\mathrm K},D(H_{\mathrm K}))$ is the  Friedrichs extension of  $(H_{\mathrm K},\FC)$
in $L^2(\Gamma,\mu)$.

\end{theorem}

\noindent {\it Proof}. By \eqref{mecke} and \eqref{hyhgfuy}, the theorem easily follows from our assumptions \eqref{igyig}, \eqref{jkg}.\quad $\square$

\subsection{Examples}\label{lkhuhy}

For each $s\in[0,1]$, we define
\begin{equation}\label{c} c(x,y,\gamma):=a(x,y)r(x,\gamma)^{s-1}r(y,\gamma)^s\One_{\{r(x,\gamma)>0,\, 
r(y,\gamma)>0\}}.
\end{equation}
Here, the function $a:X^2\to[0,\infty)$ is measurable, symmetric (i.e., $a(x,y)=a(y,x)$), bounded, and satisfies
\begin{equation}\label{uyfuyt}
\sup_{x\in X}\int_X a(x,y)\,\nu(dy)<\infty.
\end{equation}
 Assume also that there exists $\Lambda\in{\cal B}_0(X)$ such that 
\begin{equation}\label{ytdtrd}
\sup_{x\in X\setminus \Lambda} J(x,x)<\infty.\end{equation}
 
  Note that $c(x,y,\gamma)$ satisfies the balance condition \eqref{hyhgfuy}

\begin{prop}\label{kjggiuiui}
Let the coefficient $c(x,y,\gamma)$ be given by  \eqref{c}, and let conditions \eqref{uyfuyt}, \eqref{ytdtrd} hold. 
Then, for each $s\in[0,1]$,  conditions 
 \eqref{gyu} and \eqref{cond2} are satisfied, and therefore the conclusion of Theorem~\rom{\ref{8435476} } holds for the corresponding Dirichlet form.
 
 Furthermore, for $s=1$, condition
\eqref{jkg} is satisfied, and hence the conclusion of Theorem \rom{\ref{gener}}
holds for the corresponding generator $(H_{\mathrm K},D(H_{\mathrm K}))$. 
\end{prop}

\noindent {\it Proof.} Let $s\in[0,1]$.
We have, by  \eqref{mecke}, \eqref{uyfuyt}, \eqref{ytdtrd} and Proposition~\ref{gtyfy},
\begin{align*}
& \int_\Gamma\mu(d\gamma)\int_X\gamma(dx)\int_X \nu(dy)c(x,y,\gamma\setminus x)(\One_\Lambda(x)+\One_\Lambda(y))\\
&\qquad=\int_\Gamma\mu(d\gamma)\int_X\nu(dx)\int_X \nu(dy)
a(x,y)r(x,\gamma)^s r(y,\gamma)^s\\
&\qquad\quad\times \One_{\{r(x,\gamma)>0,\, 
r(y,\gamma)>0\}}(\One_\Lambda(x)+\One_\Lambda(y))\\
&\qquad\le \int_\Gamma\mu(d\gamma)\int_X\nu(dx)\int_X \nu(dy)
a(x,y)\\
&\qquad\quad\times (1\wedge J(x,x))(1\wedge J(y,y))(\One_\Lambda(x)+\One_\Lambda(y))<\infty,
\end{align*} so that condition \eqref{gyu} is satisfied.

Next, setting $u=-s$, we see that in order to show that \eqref{cond2} is satisfied, it suffices to prove that, for each $\Lambda\in{\cal B}_0(X)$,
$$\int_\Lambda\nu(dx)\int_\Lambda \gamma(dy)a(x,y)r(x,\gamma\setminus y)^s\in L^2(\mu). $$
So, by    Proposition~\ref{gtyfy}, \eqref{mecke},  \eqref{uyfuyt}, \eqref{ytdtrd}, and the boundedness of $a$,   we have:
\begin{align*}
&\int_\Gamma\mu(d\gamma)\bigg(\int_\Lambda\nu(dx)\int_\Lambda \gamma(dy)a(x,y)r(x,\gamma\setminus y)^s\bigg)^2\\
&\qquad=\int_\Gamma\mu(d\gamma)\int_\Lambda \nu(dy)r(y,\gamma)
\int_\Lambda\nu(dx_1)\int_\Lambda\nu(dx_2)a(x_1,y)a(x_2,y)
r(x_1,\gamma)^sr(x_2,\gamma)^s\\
&\qquad\quad+
\int_\Gamma\mu(d\gamma)\int_\Lambda\nu(dy_1)\int_\Lambda\nu(dy_2)\int_\Lambda\nu(dx_1)\int_\Lambda\nu(dx_2)r(y_2,\gamma)r(y_1,\gamma\cup y_2)\\
&\qquad\qquad\times a(x_1,y_1)a(x_2,y_2)r(x_1,\gamma\cup y_2)^s
r(x_2,\gamma\cup y_1)^s\\
&\qquad\le \int_\Lambda \nu(dy)J(y,y)
\int_\Lambda\nu(dx_1)\int_\Lambda\nu(dx_2)a(x_1,y)a(x_2,y)
(1+J(x_1,x_1))(1+J(x_2,x_2))\\
&\qquad+\int_\Lambda\nu(dy_1)\int_\Lambda\nu(dy_2)\int_\Lambda\nu(dx_1)\int_\Lambda\nu(dx_2) a(x_1,y_1)a(x_2,y_2)\\
&\qquad\qquad\times J(y_1,y_1)J(y_2,y_2)(1+J(x_1,x_1))(1+J(x_2,x_2))<\infty.\end{align*}

Now, let $s=1$. Analogously to the above, we have:
\begin{align*}
&\int_\Gamma \mu(d\gamma)\bigg(\int_X\gamma(dx)\int_X\nu(dy)c(x,y,\gamma\setminus x)(\One_\Lambda(x)+\One_\Lambda(y))\bigg)^2\\
&\qquad =\int_\Gamma\mu(d\gamma)\int_X\nu(dx)r(x,\gamma)
\int_X\nu(dy_1)\int_X\nu(dy_2) a(x,y_1)a(x,y_2)\\
&\qquad\quad\times r(y_1,\gamma) r(y_2,\gamma)\One_{\{r(x,\gamma)>0,\, 
r(y_1,\gamma)>0,\, r(y_2,\gamma)>0\}}
(\One_\Lambda(x)+\One_\Lambda(y_1))(\One_\Lambda(x)+\One_\Lambda(y_2))\\
&\qquad \quad+\int_\Gamma\mu(d\gamma)\int_X\nu(dx_1)\int_X\nu(dx_2 )r(x_2,\gamma)r(x_1,\gamma\cup x_2)\\
&\qquad\qquad\times
\int_X\nu(dy_1)\int_X\nu(dy_2) a(x_1,y_1)a(x_2,y_2) r(y_1,\gamma\cup x_2) r(y_2,\gamma\cup x_1)\\
&\qquad\qquad\times\One_{\{r(x_1,\gamma\cup x_2)>0,\,
r(x_2,\gamma\cup x_1)>0,\, r(y_1,\gamma\cup x_2)>0,\, r(y_2,\gamma\cup x_1)>0\}}\\
&\qquad\qquad\times(\One_\Lambda(x_1)+\One_\Lambda(y_1))(\One_\Lambda(x_2)+\One_\Lambda(y_2))\\
&\qquad\le \int_X\nu(dx)
\int_X\nu(dy_1)\int_X\nu(dy_2) a(x,y_1)a(x,y_2)\\
&\qquad\quad\times 
J(y_1,y_1)J(y_2,y_2)
(\One_\Lambda(x)+\One_\Lambda(y_1))(\One_\Lambda(x)+\One_\Lambda(y_2))\\ &\qquad\quad+\int_X\nu(dx_1)\int_X\nu(dx_2 )\int_X\nu(dy_1)\int_X\nu(dy_2) a(x_1,y_1)a(x_2,y_2)\\
&\qquad\qquad\times J(y_1,y_1)J(y_2,y_2)(\One_\Lambda(x_1)+\One_\Lambda(y_1))(\One_\Lambda(x_2)+\One_\Lambda(y_2))<\infty.\quad\square
\end{align*}

Next, for each $s\in[0,1]$, we define
\begin{equation}\label{jkgbjhgv}
d(x,\gamma):=r(x,\gamma)^{s-1}\One_{\{r(x,\gamma)>0\}},\end{equation}
so that $$ b(x,\gamma):=r(x,\gamma)^{s}\One_{\{r(x,\gamma)>0\}}. $$

Analogously to Proposition \ref{kjggiuiui}, we get

\begin{prop} Let the coefficient $d(x,\gamma)$ be given by \eqref{jkgbjhgv}. Then, for each $s\in[0,1]$, condition \eqref{gyiu}, is satisfied, and hence the conclusion of Theorem~\rom{\ref{8435476} } holds for the corresponding Dirichlet form. 

Furthermore, for $s=1$, condition \eqref{igyig} is satisfied, 
and hence the conclusion of Theorem \rom{\ref{gener}}
holds for the corresponding generator $(H_{\mathrm G},D(H_{\mathrm G}))$. 
\end{prop} 

We finally note that all our assumptions are trivially satisfied in the case of bounded coefficients $c(x,y,\gamma)$ and $d(x,\gamma)$, $b(x,\gamma)$, respectively.

\begin{center}
{\bf Acknowledgements}\end{center}

We are grateful to Yuri Kondratiev for many fruitful discussions. We would like to thank 
Dmitri Finkelstein,  Tobias Kuna, and Oleksandr Kutoviy for useful comments.
 N.T. gratefully acnowledges the financial support by the DFG-Graduiertenkolleg 
"Stoch\-astics and Real World Models". E.L. gratefully acknowledges the financial support of the DFG 
through SFB 701, Bielefeld University.

\end{document}